\RequirePackage{fix-cm}
\documentclass[smallcondensed]{svjour3}     
\smartqed  
\usepackage{graphicx}
\usepackage{comment}

\usepackage{algorithm}
\usepackage{algpseudocode}
\usepackage{amsfonts}
\usepackage[fleqn]{amsmath}
\usepackage{amssymb}
\usepackage{array}
\usepackage{booktabs}
\usepackage{cancel}
\usepackage{color}
\usepackage{comment}
\usepackage{enumerate}
\usepackage{epsfig}
\usepackage{epstopdf}
\usepackage{fancybox}
\usepackage{float}
\usepackage{footnote}
\usepackage[margin=1in]{geometry}
\usepackage{graphics}
\usepackage{graphicx} 
\usepackage{latexsym}
\usepackage{lineno}
\usepackage{lscape}
\usepackage{multirow}
\usepackage[comma, numbers]{natbib} 

\usepackage{pgfplotstable}
\pgfplotsset{compat=1.14}

\usepackage{psfrag}
\usepackage{rotating}
\usepackage{tabularx}
\usepackage[flushleft]{threeparttable}
\usepackage{url}

\algtext*{EndIf}
\algtext*{EndFor}
\algtext*{EndFunction}
\algtext*{EndWhile}

\newcolumntype{P}[1]{>{\centering\arraybackslash}p{#1}}

\newcommand{\lp}{\left(}
\newcommand{\rp}{\right)}

\newcommand{\update}[1]{\textcolor{black}{#1}}

\renewcommand{\gets}{\leftarrow}

\newtheorem{defn}{Definition}

\begin{document}
	\title{Robust flight schedules with stochastic programming
	}
	
	
	\author{Sujeevraja Sanjeevi$^\text{a}$ \and Saravanan Venkatachalam$^\text{b}$$^\dagger$}
	
	
    \institute{
        $^\text{a}$ \email{sujeev.sanjeevi@gmail.com}\\
        $^\text{b}$Department of Industrial and Systems Engineering, Wayne State University, Detroit, MI, USA.\\
        $^\dagger$Corresponding author, \email{saravanan.v@wayne.edu}
}
	
	\date{Received: date / Accepted: date}
	
	\maketitle
	
	\begin{abstract}
		Limiting flight delays during operations has become a critical research topic in recent years due to their prohibitive impact on airlines, airports, and passengers. A popular strategy for addressing this problem considers the uncertainty of day-of-operations delays and adjusts flight schedules to accommodate them in the planning stage. In this work, we present a stochastic programming model to account for uncertain future delays by adding buffers to flight turnaround times in a controlled manner. Specifically, our model adds slack to flight connection times with the objective of minimizing the expected value of the total propagated flight delay in a schedule. We also present a parallel solution framework that integrates an outer approximation decomposition method and column generation. Further, we demonstrate the scalability of our approach and its effectiveness in reducing delays with an extensive simulation study of five different flight networks using real-world data.
		
	\end{abstract}
	
	\section{Introduction}
	\label{intro}
	The Bureau of Transportation Statistics reports that between October $2018$ and October $2019$, delays caused by late aircraft arrivals amounted to 40+ million minutes, which is $39.47\%$ of the total delays experienced by the flights of reporting carriers \cite{bts}. This highlights that operational delays are a significant problem on both an absolute and a relative basis even today, with propagated delays being the biggest offender. Propagated delays occur when the arriving flight for a connection is delayed and causes a departure delay for the onward flight, kicking off a chain reaction of delays on the aircraft's route. Such propagation is primarily due to the creation of ``tight'' schedules with very limited buffers for connection times. Such schedules are created to maximize utilization of assets such as equipment and crew \cite{klabjan2001robust}. This leaves no room for the schedule to absorb fluctuations in flight arrivals and departures, resulting in significant delays and costs.
	
	The idea of making an airline schedule robust seeks to counteract this problem by adjusting the schedule to better absorb time fluctuations in aircraft arrivals and departures during operations. As robustness-based decisions need to be made much earlier than actual operational delays are known, it is necessary to consider the stochasticity of such delays.  The downside of this approach is a reduction in resource utilization and an increase in planned operational costs. This creates the need for solution strategies that can balance planning and operational costs. Optimization-based approaches, which are inherently equipped with mechanisms for such balancing acts, are therefore a great fit for this problem.
	
	Schedule robustness has been tackled in the literature from several perspectives. A two-stage stochastic programming model is proposed in \cite{yen2006stochastic}, where crew assignments are made in the first-stage and swap opportunities are anticipated in the second-stage. Another two-stage stochastic programming model is presented in \cite{froyland2013recoverable}, where the first-stage is a tail assignment problem and the second-stage is a schedule recovery problem. This model uses penalties to minimize changes between the planning and recovery solutions. A mixed integer program (MIP) with stochastic data to minimize expected propagated delays is presented in \cite{lan2006planning}. The study in \cite{marla2018robust} compares the performance of chance-constrained programming, robust optimization, and stochastic optimization approaches using a solution space similar to the one in the model presented in \cite{lan2006planning}. Methodologies to solve integrated aircraft routing and crew pairing problems to reduce uncertain propagated delays are considered in \cite{yen2006stochastic,dunbar2012robust,dunbar2014integrated}. More recently, the robust optimization approach presented in \cite{yan2016robust} uses column and row generation to solve a routing problem with delays coming from a bounded uncertainty set by minimizing worst-case propagated delay costs. An alternate perspective in \cite{ahmadbeygi2010decreasing,chiraphadhanakul2013robust} retains a given planned routing but re-times flights in order to add time buffers or ``connection slacks'' to flight connections that are likely to be missed. Other related work can also be found in \cite{arikan2013building,kang2004degradable,rosenberger2004robust,shebalov2006robust,talluri1996swapping,weide2010iterative}.
	
	\update{To motivate our research, we present some concerns we observed with scheduled robustness models proposed so far. First, there is no clear differentiation between the cost of rescheduling flights a few weeks before the day-of-operations versus delaying them a few hours before departure. This difference can be significant in practice. Second, the stochastic programming approaches proposed in literature use very complex first-stage models with a wide variety of first-stage decisions. This may be undesirable, as each adjustment of a schedule can affect other operational considerations such as staff scheduling, maintenance scheduling, crew and passenger connectivity, among others. Also, there is no clarity on how to reduce the scope of such models while still generating useful results for scheduling practitioners. Computationally, the size and complexity of first-stage models proposed in literature makes it difficult to scale them and use them for real-world airline schedules.}
	
	\update{In this research, we seek to fill the aforementioned gaps in literature. Our main contributions are (i) a two-stage stochastic programming model that re-times flights of a given schedule in a controlled manner while minimizing the sum of first-stage rescheduling costs and expected cost of propagated delays on the day of operations; (ii) a parallel decomposition framework based on the L-shaped method \cite{van1969shaped} that uses column generation to solve recourse second-stage models; (iii) extensive computational study using disruptions caused by randomly generated flight delays that show a significant reduction in propagated delays in schedules adjusted by the proposed model; and (iv) recommendations and insights  that can boost the performance of decomposition techniques like Benders decomposition and column generation for flight schedule models. The proposed model and solution framework allow to solve much larger instances than those solved so far in literature. For example, one of the networks we consider has $324$ flights and $71$ aircraft, much larger in size than networks used in recent works like \cite{froyland2013recoverable}, \cite{yan2016robust}. Furthermore, we use a dynamic delay approach similar to \cite{yan2016robust} to solve our recourse problems. This approach uses the least required delay on each flight while building paths. This eliminates the need for discrete delay copies which can generate unnecessary flight delays due to discretization and cause significant run time increases (see Figure 7 in \cite{froyland2013recoverable}). The path-based recourse formulation in our model can be easily extended to incorporate requirements from other operational domains of airlines. This includes hard constraints like minimum crew/passenger connection times and soft requirements like the loss of passenger goodwill that can be incorporated into path costs.}

	The remainder of this paper is organized as follows. In Section \ref{sec:models}, we present a two-stage stochastic programming formulation to minimize the expected value of propagated delays, along with a simpler mixed integer programming formulation based on sample mean values of primary delays. In Section \ref{sec:sol}, we describe a column-generation procedure for recourse problems and the L-shaped algorithm for the complete two-stage problem. In Section \ref{sec:comp}, we report the results of extensive computational studies that highlight the qualitative and quantitative benefits of our approach. In Section \ref{sec:conclusions}, we conclude the article with a summary and  discussion of future research directions.
	
	\section{Stochastic Delay Models}\label{sec:models}
	
	In this section, we present our two-stage stochastic programming formulation of the delay mitigation problem. We also present an alternate approach that we use to benchmark our computational results. The latter approach is based on an MIP model that uses the mean values of individual flight delays. We begin by introducing the required notation.
	
	Given a valid flight schedule, we model it as a connection network on a directed acyclic graph $G = (F, A)$ in which the set of nodes $F$ represent flights and the arcs $A$ represent flight connections. A connection $(i,j)$ is valid if and only if (i) the incoming arrival and outgoing departure airports match, and (ii) the connection slack $s_{ij}$, defined as the difference between the departure time of the outgoing flight $j$ and the arrival time plus the turnaround time of the incoming flight $i$, is non-negative. The set $A$ contains only valid connections.
	
	Our modeling of uncertain flight delays is similar to that in \cite{lan2006planning,dunbar2012robust,yan2016robust}. A flight can experience \textit{primary delays} that are independent of routing and rescheduling, and \textit{propagated delays} that are caused by upstream flights on that flight's route. Let $\omega$ be a random variable representing a delay scenario, and let $\Omega$ be a finite set of delay scenarios. Let $pd_f^\omega$ be the realized non-negative integer-valued primary delay in minutes experienced by flight $f \in F$ in scenario $\omega \in \Omega$. Let $R^\omega$ be the set of possible routes in scenario $\omega$. For any route $r \in R^\omega$ and connection $(i,j)$ in $r$, the parameter $d_{rj}$ representing the delay propagated to the outgoing flight $j$ by the connection is defined as:
	\begin{align}
		d_{rj} = max(0, d_{ri} + pd_i^\omega - s_{ij}). \label{eq:propagatedDelay}
	\end{align}
	
	\subsection{Two-stage model}\label{subsec:two-stage}
	Let $x_f \geq 0$ be an integer decision variable representing the number of minutes by which flight $f \in F$ needs to be rescheduled, and let $c_f, f \in F$, be the per-minute reschedule cost. The formulation of the two-stage model (TSM) can then be stated as:
	\begin{align}
		\nonumber (TSM) \quad & \text{Minimize } && \sum_{f \in F} c_f x_f + \mathbb{E}_\Omega[\phi(x, \tilde{\omega})] && \\
		\label{eq:origRoute} & \text{s.t. } && x_i \leq s_{ij} + x_j, && (i,j) \in A^{orig}, \\
		\label{eq:rescheduleBudget} & && \sum_{f \in F} x_f  \leq B, && \\
		\label{eq:firstvars} & && x_f \in \mathbb{Z} \cap [0,l], && f \in F.
	\end{align}
	
	\noindent The objective of this model is to minimize the sum of the total reschedule cost and the expected flight delay costs. Constraints \eqref{eq:origRoute} protect the time connectivity for all connections in the original routing $A^{orig} \subseteq A$. Constraints \eqref{eq:rescheduleBudget} provide a control factor in the form of a time budget $B$ that limits the total reschedule time. We also limit the $x_f$ values with a fixed bound $l$ to prevent exorbitant reschedules of individual flights. Given a reschedule $x$ and the scenario probabilities $p_\omega$, $\omega \in \Omega$, the expected value $\mathbb{E}_\Omega[\phi(x, \omega)] = \sum_{\omega \in \Omega} p_\omega \phi(x, \omega)$ can be computed by solving the following set partitioning model for each scenario $\omega \in \Omega$, which is the second-stage formulation for a given $x$ and scenario $\omega$:
	\begin{align}
		\phi(x, \omega) = \nonumber& \text{ Min } && \sum_{f \in F} e_f z_f^\omega && && \\
		\label{eq:onePerTail} & \text{ s.t. } && \sum_{r \in R^\omega} a_{rt} y_r  = 1, && t \in T,\\
		\label{eq:onePerFlight} & && \sum_{r \in R^\omega} b_{rf} y_r  = 1, && f \in F,\\
		\label{eq:delayLink} & && \sum_{r \in R^\omega} b_{rf} d_{rf} y_r - x_f  \leq z_f^\omega, && f \in F,\\
		\label{eq:ssvars}& && z_f^\omega \geq 0, \text{ } f \in F, y_r \in \{0, 1\}, \ r \in R^\omega. &&
	\end{align}
	
	\noindent The second-stage model minimizes the propagated delay costs incurred in scenario $\omega \in \Omega$ computed as per-minute costs $e_f$ for each flight $f$. It uses two sets of decision variables: continuous variables $z_f^\omega$ that represent the excess delay propagated to each flight $f \in F$ and binary variables $y_r$ that take the value $1$ to indicate the selection of the route $r \in R^\omega$. The parameters $a_{rt}$ and $b_{rf}$ are binary  and respectively indicate whether route $r$ is for the tail $t$ and whether it contains flight $f$. Constraints \eqref{eq:onePerTail} and \eqref{eq:onePerFlight} enforce the assignment of one route per aircraft and one route per flight. Constraints \eqref{eq:delayLink} are linking constraints that capture the excess propagated delay that has not been accounted for by the first-stage rescheduling.

	Next, we present an MIP formulation that reschedules flights based on the average values of the primary delays. This model is used in the comparative studies presented in the computational results section.
	
	\subsection{Mean delay model}
	Let $\bar{\omega}$ be the scenario in which each flight experiences the mean primary delay across all scenarios in $\Omega$, i.e., $d_f^{\bar{\omega}} = \sum_{\omega \in \Omega} p_\omega d_f^\omega$ for $f \in F$. The mean delay model aims to reschedule flights to accommodate the average delay scenario $\bar{\omega}$ without changing the original routing. To simplify the notation, we set $d_f^{\bar{\omega}}$ to be the delay propagated to flight $f$ in scenario $\bar{\omega}$ in the original schedule. The mean delay model can be stated as follows:
	\begin{align*}
		& \text{Minimize } && \sum_{f \in F} \lp c_f x_f + e_f z_f^{\bar{\omega}} \rp && && \\
		& \text{s.t. } && x_i \leq s_{ij} + x_j, && (i,j) \in A^{orig},\\
		& && \sum_{f \in F} x_f  \leq B, && \\
		& && d_f^{\bar{\omega}} - x_f  \leq z_f^{\bar{\omega}}, && f \in F,\\
		& && z_f^{\bar{\omega}} \geq 0, x_f \in \mathbb{Z} \cap [0,l], && f \in F.
	\end{align*}
	
	\noindent The objective function minimizes the total reschedule and delay costs, with the latter carrying a higher penalty. The first two sets of constraints are the first-stage constraints \eqref{eq:origRoute} and \eqref{eq:rescheduleBudget}. The third set of constraints is obtained from \eqref{eq:delayLink} by selecting only the original route for each aircraft.
	
	\section{Solution approach}\label{sec:sol}
	In this section, we present our solution framework that uses the $L$-shaped method in \cite{van1969shaped} to solve the TSM. We first present details about how we solve the recourse problems of the TSM.
	
	\subsection{Column-generation framework}
	Solving the TSM using the $L$-shaped method requires computing $\phi_{LP}(\bar{x}, \omega)$, the solutions to linear programming (LP) relaxations of the recourse models for any fixed first-stage solution $\bar{x}$. For a given \update{scenario} $\omega$, we use a column-generation approach to generate the required routes. We iterate between solving a version of the recourse problem restricted to a subset of routes $\tilde{R} \subseteq R^\omega$ and solving a pricing problem to find new routes that can improve the solution. \update{Optimality of the linear program} can be declared when no such route can be found. For ease of exposition, we state here the dual formulation of the recourse problem in full. Let $\mu_t$ and $\nu_f$ be unbounded dual variables for the coverage constraints \eqref{eq:onePerTail} and \eqref{eq:onePerFlight} for a scenario $\omega$. Given a first-stage solution $\bar{x}$, we write the constraints \eqref{eq:delayLink} as
	\[ z_f^\omega - \sum_{r \in R^\omega} b_{rf} d_{rf} y_r \geq -\bar{x}_f,\ f \in F, \]
	
	\noindent and we let $\pi_f$ be the non-negative dual variables for these constraints. Let $a(r) \in T$ be the aircraft for which the route $r \in R^\omega$ was generated. Using this notation, the dual formulation can be written as:
	\begin{align}
		\nonumber
		& \text{ Maximize } && \sum_{t \in T} \mu_t + \sum_{f \in F} \lp \nu_f - \bar{x}_f \pi_f \rp && \\
		\label{eq:dualRoute}
		& \text{ s.t. } && \mu_{a(r)} + \sum_{f \in F} b_{rf} \lp\nu_f - d_{rf} \pi_f \rp \leq 0, && r \in R^\omega,\\
		\nonumber
		& && \mu_t \text{ free}, && t \in T,\\
		\nonumber
		& && \nu_f \text{ free, } 0 \leq \pi_f \leq e_f, && f \in F.
	\end{align}
	
	Our column-generation procedure begins by solving the LP relaxation of the recourse problem with a subset $\tilde{R}$ of routes. One way to initialize $\tilde{R}$ is with routes of the original schedule that have delays propagated sufficiently enough to protect minimum turnaround times. With the dual solution of this restricted problem, a pricing problem is solved to find columns with the least reduced cost $rc_r$, where
	\begin{equation}
		rc_r = \sum_{f \in F} b_{rf} \lp d_{rf} \pi_f - \nu_f \rp - \mu_{a(r)}.
		\label{eq:reducedCost}
	\end{equation}
	
	\noindent The dual formulation provides some intuition for $rc_r$; we want routes that violate the constraints \eqref{eq:dualRoute}. Once such a route is found, it is added to $\tilde{R}$ and we repeat the above steps. If no such route can be found, optimality can be declared. As there are potentially a large number of pricing problems to be solved, it is critical to determine the useful routes quickly. Next, we present our version of the labeling algorithm, an extension of the algorithm presented in \cite{dunbar2012robust,yan2016robust}, which we use to solve this problem.
	
	\subsection{Pricing problem} \label{subsec:pricing-problem}
	We solve the pricing problem by searching for routes in the graph $G$ with negative values for the reduced cost as defined in \eqref{eq:reducedCost}. As we assume that the original schedule is already available, the airports from which each aircraft should depart at the beginning of the schedule and at which it should arrive at the end of the schedule are fixed. To reflect this, we introduce separate source and sink nodes for each aircraft and separately search for candidate routes for each aircraft. This approach is quite practical, as it can easily be extended to consider aircraft-specific business constraints during route generation. Each aircraft's source node connects only to flights departing from the aircraft's initial departure airport. Similar restrictions apply to sink nodes based on final arrival airports.
	
	To search for candidate routes, we use a label-setting algorithm similar to the one proposed in \cite{dunbar2012robust,yan2016robust}. This algorithm relies on building \textit{labels} that represent partial routes and extending them along valid flight connections given by $A$ to generate full routes from the source to the sink. The combinatorial explosion in the number of routes is controlled using the notion of \textit{dominance} between labels. More formally, each label $l$ denotes a partial path stored in a tuple $(f_l, pred_l, red_l, prop_l)$, where $f_l \in F$ is the last flight on the path, $pred_l$ is the label from which $l$ was extended, $red_l$ is the reduced cost accumulated so far, and $prop_l$ is the delay propagated to $f_l$ on the partial route corresponding to $l$. Note that $pred_l$ is empty for labels at source nodes. When a label $u$ is extended with a connection $(f_u, f') \in A$, the algorithm generates a new label $v = (f',u, red_v,prop_v)$ in which $red_v$ and $prop_v$ are updated using \eqref{eq:propagatedDelay} and \eqref{eq:reducedCost}, respectively. Once a label is extended to the sink node, the route that it corresponds to becomes a full route and can be obtained by traversing backward along the chain of predecessors.
	
	\begin{defn} \label{def:dominance}
		(Label dominance condition) Let $u$ and $v$ be two labels with $f_u=f_v$. The label $u$ dominates $v$ if (i) $red_u \leq red_v$, (ii) $prop_u \leq prop_v$, and at least one of the inequalities is strict.
	\end{defn}
	
	Given two labels $u$ and $v$, if we know that any feasible extension of $v$ is also feasible for $u$, any route that can be generated by successively extending $v$ to the sink can also be generated by $u$, meaning that we can safely ignore $v$. This was proved in Lemma $1$ in \cite{yan2016robust}. For clarity, we restate the lemma here using the notation of the present article:
	\begin{lemma} \label{lem:extention}
		Let $u$ and $v$ be labels such that $u$ dominates $v$. If $u'$ and $v'$ are labels obtained by extending $u$ and $v$ with a connection $(f_u,f') \in A$, then $u'$ dominates $v'$.
	\end{lemma}
	
	Lemma \ref{lem:extention} allows us to store and extend only non-dominated labels at each node and thus implicitly remove large numbers of candidate paths from consideration. We have observed that the label-setting algorithm in \cite{yan2016robust} provides at most one negative reduced-cost route in each iteration. As any route with a negative reduced cost is likely to improve the recourse solution, we enhance the algorithm by considering three possible alternatives for generating multiple negative reduced-cost columns:
	
	\begin{enumerate}[(i)]
		\item \textit{All paths}: Store and return all negative reduced-cost paths.
		\item \textit{Best paths}: Store all negative reduced-cost paths, but return only the $N$ most negative reduced-cost paths.
		\item \textit{First paths}: Stop the search as soon as $N$ negative reduced-cost paths are found, and return them.
	\end{enumerate}
	
	We found that all three strategies produce a significant speedup over generating a single path per pricing problem. Among the three, the ``first paths'' strategy gave us the best runtime with $N$=10. We present a more detailed comparative study of these strategies in the computational results section. We present the label-setting algorithm of \cite{dunbar2012robust,yan2016robust} with our enhancements below, in Algorithm \ref{algo:pricing}. As the original initial-departure and final-arrival airports can be different for each aircraft, the algorithm is used to separately generate routes for each aircraft. The input includes augmented sets of nodes $F'$ and arcs $A'$; $F'= F \cup \{so,si\}$, where $so$ and $si$ are dummy source and sink nodes, respectively, and $A'$ contains all eligible connections in $A$, connections from $so$ to every valid first flight in $F$, and connections from every valid last flight to $si$ for the selected aircraft. The output of the algorithm is a set of negative reduced-cost columns for the selected aircraft.
	
	\begin{algorithm}[!htbp]
		\caption{Label-setting algorithm}\label{algo:pricing}
		\begin{algorithmic}
			\Function{GenerateColumns}{$F'$, $A'$, $so$, $si$}
			\State $M_f \gets \emptyset$, $f \in F'$. \Comment{Processed labels container}
			\State $I_f \gets \emptyset$, $f \in F'$. \Comment{All labels container}
			\State $I_{so} \gets \{(so, \varnothing, -\mu_{a(r)}^\omega, 0)\}$. \Comment{Source label creation}
			\While{$\bigcup_{i \in F'}(I_i\setminus M_i) \neq \emptyset$ and \Call{ShouldStop}{$I_{so}$} $\neq true$}
			\State Choose $i \in F'$ and a label $l \in I_i \setminus M_i$ with a minimal reduced cost.
			\For{$(i,j) \in A'$}
			\State $l' \gets$ \Call{Extend}{l, j}.
			\If{$l'$ is not dominated by any label in $I_j$}
			\State $I_j \gets I_j \cup \{r(i),j\}$. 
			\EndIf
			\If{$j = si$ and \Call{ShouldStop}{$I_{so}$} $= true$}
			\State \textbf{break}. \Comment{Stop processing labels}
			\EndIf
			\EndFor
			\State $M_i \gets M_i \cup \{l\}$.
			\EndWhile
			\State \Return \Call{BuildColumnsFromLabels}{$I_{so}$}
			\EndFunction
		\end{algorithmic}
	\end{algorithm}
	
	Algorithm \ref{algo:pricing} initializes a single label at the source node as $(so, \varnothing, -\mu_{a(r)}^\omega, 0)$, without a predecessor. Given a label $l = (i, pred_l, red_l, prop_l)$ and a connection $(i,j)$, the \textsc{Extend} procedure creates a new label $l'$ at node $j$ by updating $prop_{l'}$ using \eqref{eq:propagatedDelay} and the reduced cost $red_{l'} = red_l + d_j \pi_j - \nu_j$, as obtained from \eqref{eq:reducedCost}. Labels become complete when they are extended to $si$. The implementation of \textsc{ShouldStop} depends on the column-generation strategy that is used. It always returns $false$ for the all-paths and best-paths strategies. For the first-paths strategy, it returns $true$ if the number of negative reduced-cost labels at $si$ have exceeded $N$, and $false$ otherwise. When the while loop ends, the \textsc{BuildColumnsFromLabels} procedure builds columns using negative reduced-cost labels at $si$. It returns all columns for the all-paths strategy, and the $N$ most negative reduced-cost columns for the other two strategies. The LP solution to the recourse problem is optimal if Algorithm \ref{algo:pricing} returns an empty set.
	
	\subsection{Solution framework for the TSM}\label{subsec:solution}
	Now that we have established the machinery to solve recourse models, we are ready to present the $L$-shaped method to solve the TSM. The method has two variants: a single-cut and a multi-cut version. We present the multi-cut method here and show later in this section how it can be modified to obtain the single-cut method. The multi-cut $L$-shaped method works with the following approximation of the TSM:
	\begin{align*}
		(MP)\quad & \text{Minimize } && \sum_{f \in F} c_f x_f + \sum_{\omega \in \Omega} \eta^\omega && \\
		& \text{s.t. } &&  \eqref{eq:origRoute} - \eqref{eq:firstvars}, \\
		&             && \eta^\omega \text{ free, } \omega \in \Omega.
	\end{align*}
	
	\noindent We refer to this version of the formulation as the ``master problem'' (MP). Our solution procedure iterates between solving the MP and the recourse LP problems. Solutions to the latter can provide \textit{optimality cuts} that bound $\eta$ from below or \textit{feasibility cuts} generated from infeasible recourse problems. As we can always get a feasible solution for any delay scenario by propagating delays along the original routing, our recourse problems are always feasible. So we only need to consider optimality cuts. To describe these cuts, we introduce the following additional notation for each scenario $\omega \in \Omega$:
	\begin{align*}
		\alpha^\omega &= p_\omega \lp \sum_{t \in T} \mu_t + \sum_{f \in F} \nu_f \rp, \text{  and  } \beta_f^\omega = p_\omega \pi_f,\ f \in F.\\
	\end{align*}
	
	\begin{algorithm}[!htbp]
		\caption{Multi-cut $L$-shaped method for the SM}
		\label{algo:l-shaped}
		\begin{algorithmic}
			\State Solve the MP without $\eta^\omega$ variables to get an initial solution $x^0$.
			\State Add $\eta^\omega$ variables to the MP.
			\State Set $UB \gets$ $\infty$, $LB \gets -\infty$, $k \gets 0$, $x^* \gets x^0$.
			\While{$UB - LB > \epsilon$ and $k \leq$ \texttt{MaxNumIterations}}
			\For{each scenario $\omega \in \Omega$}
			\State Find $\phi_{LP}(x^k,\omega)$ using column generation.
			\State Compute $\beta^\omega,\alpha^\omega$ using optimal dual values.
			\State Add cut $\eta^\omega \geq \alpha^\omega - \sum_{f \in F} \beta^\omega_f x_f$ to the MP.
			\EndFor
			\State Set $UB \gets min \lp UB, \sum_{f \in F} c_f x^k_f + \sum_{\omega \in \Omega} \phi_{LP}(x^k,\omega) \rp$.
			\If{$UB$ changed}
			\State Update incumbent solution $x^* \gets x^k$.
			\EndIf
			\State Solve the updated MP to get the objective value $obj_k$.
			\State Set $LB \gets max(LB, obj_k)$, $k \gets k+1$.
			\EndWhile
			\Return $x^*$.
		\end{algorithmic}
	\end{algorithm}
	
	Using this notation, the multi-cut procedure is presented in Algorithm \ref{algo:l-shaped}. We found that $x^0_f = 0, f \in F$ is a reasonable starting solution. The parameter \texttt{MaxNumIterations} provides a practical way to limit the algorithm's runtime. To convert the algorithm into the single-cut $L$-shaped method, we use a single variable $\eta$ in the $MP$ and add only the single cut \eqref{eq:singleCut} that is computed using the optimal dual values of all recourse problems in each iteration:\begin{equation}
		\eta \geq \sum_{\omega \in \Omega} \alpha^\omega - \sum_{f \in F} \lp \sum_{\omega \in \Omega} \beta_f^\omega \rp x_f.
		\label{eq:singleCut}
	\end{equation}
	
	We note here that the Benders cuts are valid only when the binary restrictions of the second-stage problems are relaxed. Making our approach exact requires embedding Algorithm \ref{algo:l-shaped} in a branch-and-bound scheme that finds integer solutions to all second-stage $y_r$ variables. However, as we found that most of the optimality gap was closed in the root node, we did not explore branching. As we shall see in Section \ref{sec:comp}, even these solutions can provide rescheduling values that significantly improve the preparedness of a schedule for uncertain delays.

	\section{Computational experiments}\label{sec:comp}
	In this section, we demonstrate the efficacy of our proposed formulation and solution approach using real-world data for five flight networks. We used Java for the implementation, with CPLEX $12.9$ as the solver. The experiments were conducted on an Intel(R) Xeon(R) CPU E5-$2640$ computer with $16$ logical cores and $80$ GB RAM. We implemented parallel processing using the thread-safe Java ``actors'' provided by the Akka actor library (available at \url{https://akka.io}). \update{All code and data used for our experiments is publicly available at \url{https://github.com/sujeevraja/stochastic-flight-scheduler}.}

	\begin{table}
		\begin{center}
			\caption{Instance details}
			\label{tab:data}
			\begin{tabular}{c c c c} 
				\hline
				\textit{Instance}   &   \textit{Number of Flights}  &  \textit{Number of aircraft} &  \textit{Number of paths}    \\
				\hline
				s1   &   210  &  41  &   48,674    \\
				s2   &   248  &  67  &   20,908    \\
				s3   &   112  &  17  &   39,242    \\
				s4   &   110  &  17  &   56,175    \\
				s5   &   80  &  13  &   190,540    \\
				s6   &   324  &  71  &   113,892    \\
				\hline
			\end{tabular}
		\end{center}
	\end{table}	
	
	\subsection{Network data and experiment setup}\label{subsec:network}
	Table \ref{tab:data} presents details about the flight networks we used. \update{Each network is based on daily schedules of two different airlines on different days in early $2017$, and is the planned schedule for a single equipment type. We avoid solving multiple equipment types together as such swaps can cause operational issues like unfilled seats or passenger spillage.} Each flight in our data has a minimum turnaround time that applies to connecting flights departing after the arrival of the flight. As the costing data for our networks is quite complex, we simplify the calculations with a first-stage reschedule cost of one per minute and a recourse delay cost of $10$ per minute for each flight. This costing serves to encode the significant increase of costs incurred by operational delays as opposed to planned reschedules. \update{The ``Number of paths'' values are the maximum number of paths that can be built during column generation. To calculate them, we build a flight network and add a dummy source and dummy sink node for each aircraft based on its original first-departure and last-arrival stations. We then add dummy source arcs to flights departing from the source node station and dummy sink arcs from flights arriving at the sink node station. The number of paths for each aircraft is recursively computed as the number of paths from the aircraft's dummy source to the aircraft's dummy sink. The total number of paths is the sum of paths of all aircraft.}
	
	We simulate primary delays by constructing $30$ randomly generated delay scenarios for each run. The scenarios are generated by varying two parameters: the distribution used for delay generation and the flights that experience primary delays. We follow the recommendation of \cite{yan2016robust} in using truncated normal, gamma, and log normal distributions for primary delays, with log normal being the default. 
	We select flights that experience primary delays using two strategies, which we call ``hub'' and ``rush''. The hub strategy selects flights from a hub, which we define as the airport with the most \update{departures} in a given schedule. The rush strategy calculates the duration between the earliest departure and the latest arrival for a schedule and selects flights departing during the first quarter of the window. This idea stems from the morning runway congestion that frequently occurs in most airports. 
	Our model limits first-stage rescheduling with two control factors, an individual limit of $l$ for each flight and a limit of $B$ minutes on the total delay. We fix $l$ to $30$ minutes in all of our runs. We make $B$ adaptive to the problem data by computing the total primary flight delay for each recourse scenario, taking the average of these values, and allowing $B$ to be a fraction of the average total primary delay. \update{Unless specified otherwise, we default to $0.5$ for $B$, $LogNormal(15,15)$ as the delay distribution, ``hub'' as the flight selection strategy, the multi-cut $L$-shaped method, the first-paths column generation strategy outlined in Section \ref{subsec:pricing-problem}, and use $30$ threads to solve $30$ second-stage problems in parallel. Solution times in all tables are reported in seconds.}  
	
	\subsection{Results and insights}\label{subsec:result}
	Our computational study contains three sets of results. The first set presents the performance metrics of our algorithm, as shown in Table \ref{tab:quality}. The \textit{Strategy} column shows the strategy we use to select flights, as explained above. We report two gaps: the percentage gap computed as $100 \times (UB-LB)/UB$ from Algorithm \ref{algo:l-shaped} in the \text{Gap} column, and the optimality gap of the solution in the \textit{Opt Gap} column. To compute the latter, we first find an upper bound \textbf{ub} by fixing the first-stage reschedule values to the solution found by Algorithm \ref{algo:l-shaped}, solving all second-stage problems without relaxing the binary restrictions, and computing the objective value as the sum of the fixed reschedule cost and the mean value of the second-stage delay costs. As the objective value of the solution found by Algorithm \ref{algo:l-shaped} is a lower bound (denoted by \textbf{lb}) for the optimal solution, we report the optimality gap as $100 \times$(\textbf{ub}-\textbf{lb})/\textbf{ub}. The columns \textit{Cuts} and \textit{Iter} report the number of Benders cuts added and the number of iterations, respectively. The main takeaways from Table $1$ are that the Benders gap is almost completely closed for most instances and that the root node closes more than $90\%$ of the optimality gap. \update{We believe that the low optimality gap occurs because of the set partitioning structure in the second-stage model in TSM. As set partitioning models are known to have a property called \textit{quasi-integrality} \cite{balas1975set,balas1972set,tahir2019integral}, their linear relaxations typically yield integer solutions in most cases.}

	\begin{table}
		\begin{center}
			\caption{Solution quality and performance}
			\label{tab:quality}
			\begin{tabular}{c c c c c c c} 
				\hline
				\textit{Strategy}  &  \textit{Instance}   &   \textit{Time}  &  \textit{Gap (\%)}  &  \textit{Opt gap (\%)}  &  \textit{Cuts}  & \textit{Iter}    \\
				\hline
				Hub  &  s1   &   78.42   &  0.35  &   3.42  &  886  &  30   \\
				&  s2   &   53.94   &  2     &   3.87  &  900  &  30    \\
				&  s3   &   15.94   &  0     &   0     &  93   &  6    \\
				&  s4   &   14.04   &  0.05  &   7.61  &  304  &  15   \\
				&  s5   &   73.16   &  0     &   6.18  &  352  &  16    \\
				&  s6   &   377.3   &  3.54  &   11.85 &  900  &  30    \\
				Rush &  s1   &   90.64   &  0.09  &   7.52  &  861  &  30   \\
				&  s2   &   71.07   &  0.5   &   7.94  &  888  &  30    \\
				&  s3   &   11.73   &  0.03  &   8.75  &  79   &  4    \\
				&  s4   &    6.37   &  0     &   0.41  &  115  &  6   \\
				&  s5   &   47.92   &  0     &   0.09  &  188  &  8    \\
				&  s6   &   144.34  &  0.04  &   1.82  &  302  &  13    \\
				\hline
			\end{tabular}
		\end{center}
	\end{table}
	
	\begin{figure}
		\centering
		\includegraphics[scale=.75]{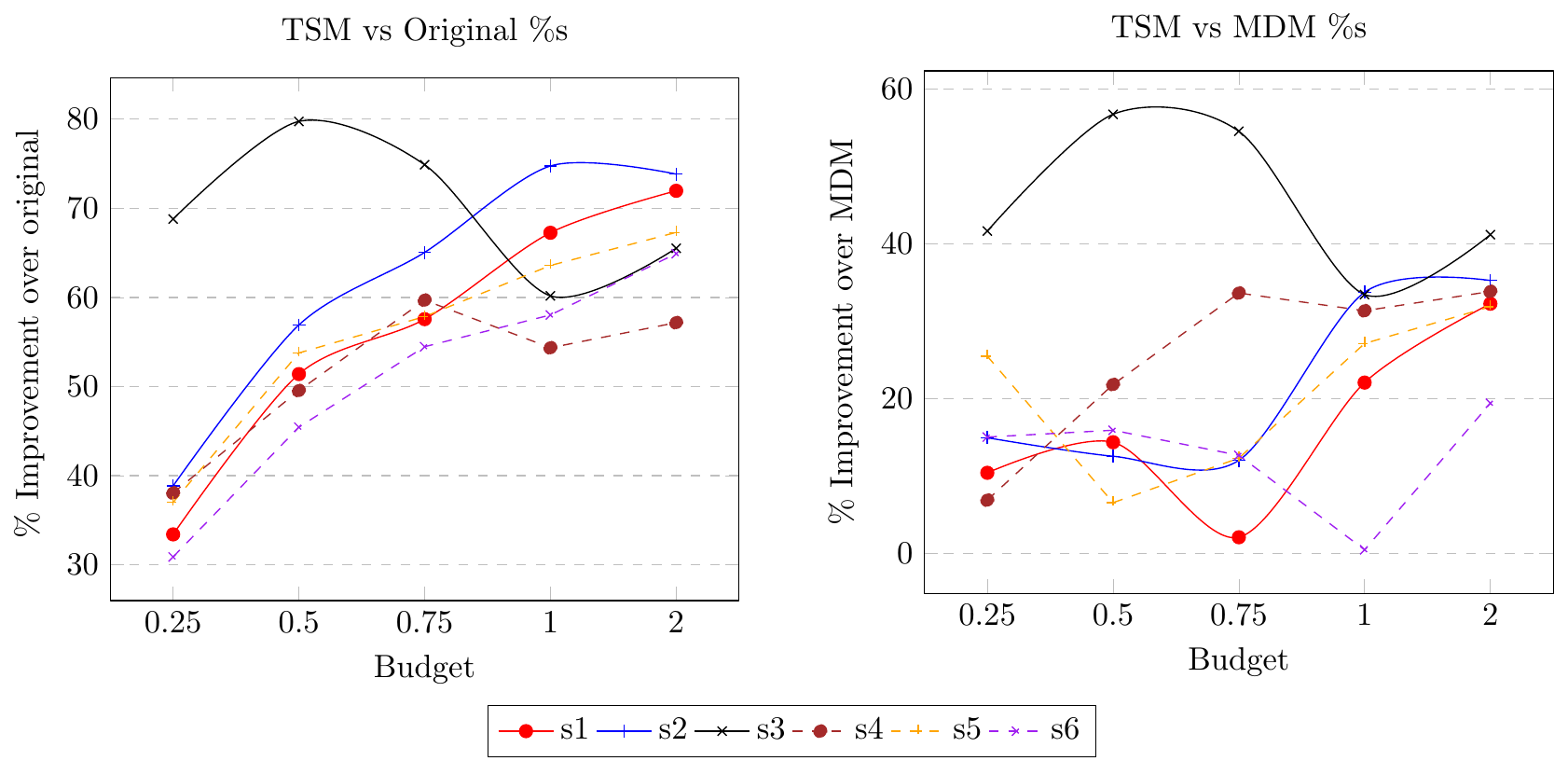}
		\caption{Illustration of performance of TSM by budget}
		\label{fig:budget}
	\end{figure}
	
	\begin{figure}
		\centering
		\includegraphics[scale=.75]{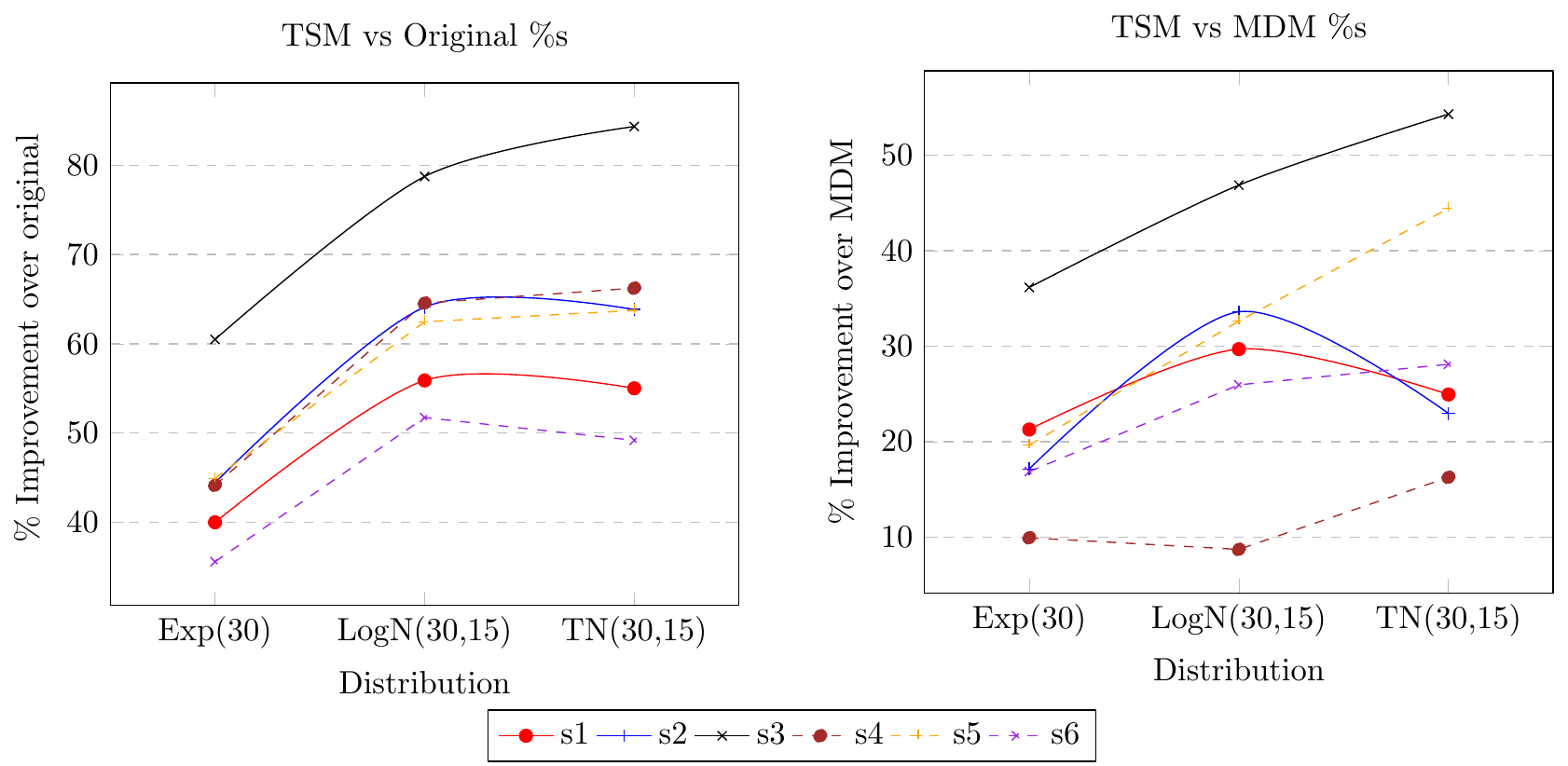}
		\caption{Illustration of performance of TSM by distribution}
		\label{fig:dist}
	\end{figure}
	
	\begin{figure}
		\centering
		\includegraphics[scale=.75]{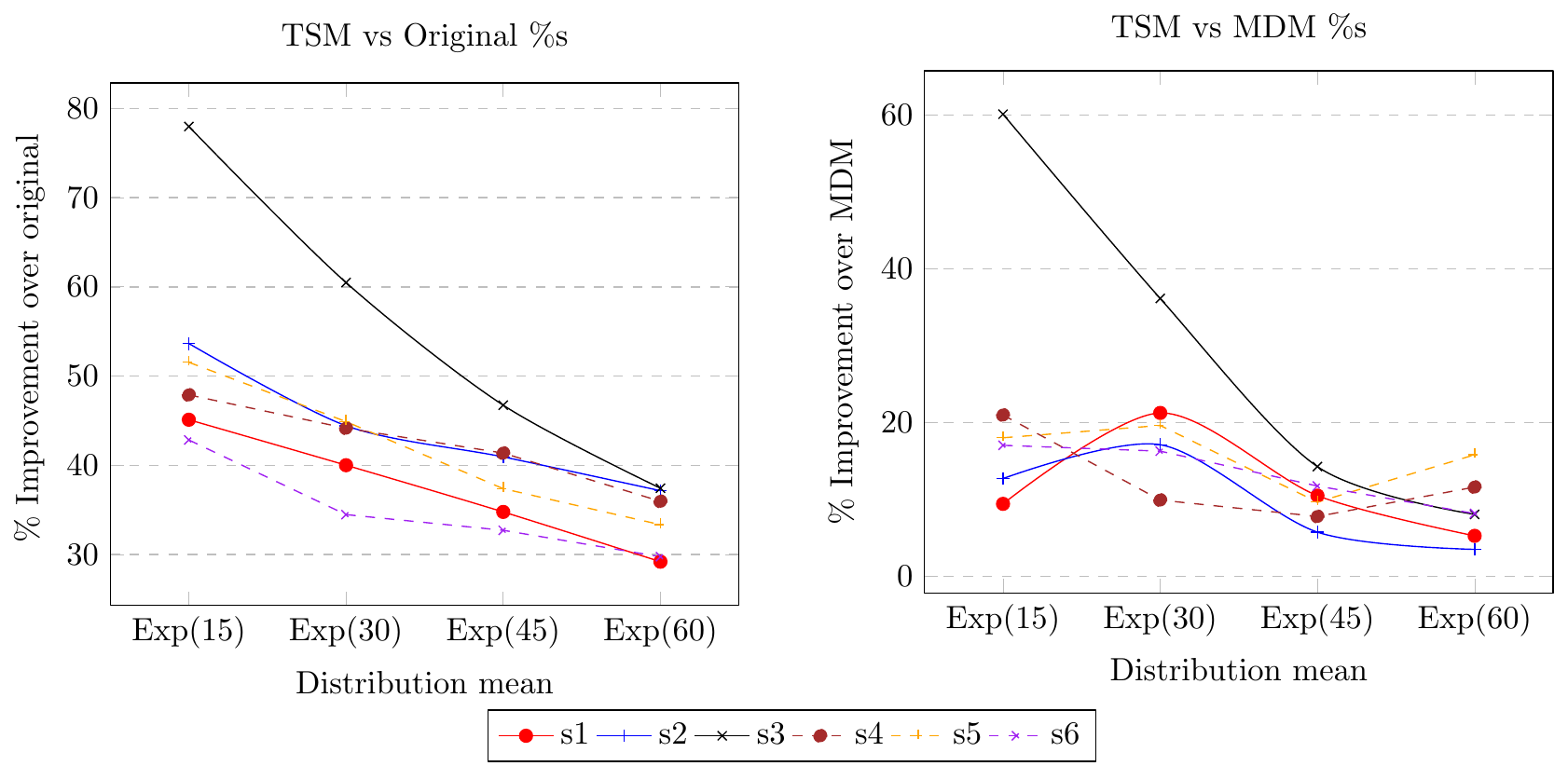}
		\caption{Illustration of performance of TSM by distribution mean}
		\label{fig:mean}
	\end{figure}
	
	\begin{figure}
		\centering
		\includegraphics[scale=.75]{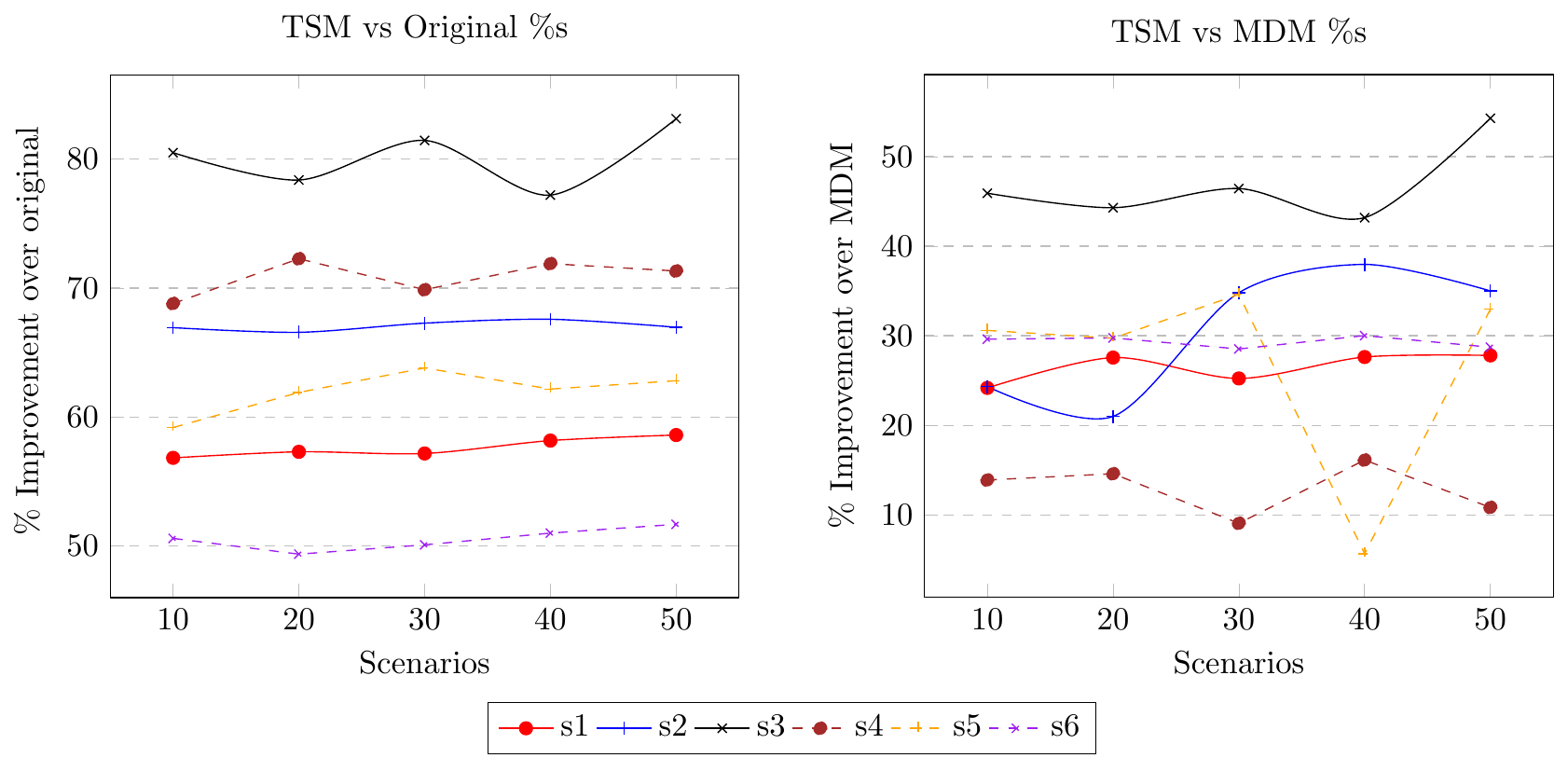}
		\caption{Illustration of performance of TSM by scenarios}
		\label{fig:scenarios}
	\end{figure}
	
	\update{For the second set of experiments, we report solution quality results in Figures \ref{fig:budget}, \ref{fig:dist}, \ref{fig:mean}, and \ref{fig:scenarios}. Numbers used for these figures are available in the Appendix in Tables \ref{tab:budget}, \ref{tab:distribution}, \ref{tab:mean} and \ref{tab:scenario} respectively.} 
	\update{In these experiments}, we first randomly generate $30$ delay scenarios and use this data to solve the two-stage and mean delay models. The same scenarios are used for both models for a fair comparison. Next, we generate a new set of $100$ random delay scenarios different from those used for solving. For each new scenario, we compute the total propagated delay incurred by three variants of the original schedule: (i) no adjustments, (ii) adjustments based on the reschedule solution of the \update{mean delay model}, and (iii) adjustments based on \update{the reschedule solution of the TSM}. By ``adjustment'', we mean that the departure time of a flight is changed based on its corresponding reschedule value. The propagated delay for any scenario\update{, measured in minutes,} is found by solving the integer-valued recourse model to optimality. We then take the average value of the total propagated delay of the $100$ scenarios as a comparison metric for the three approaches. \update{Each set of figures \ref{fig:budget}, \ref{fig:dist}, \ref{fig:mean}, and \ref{fig:scenarios} has two charts that measure the relative reduction in the average value of total propagated delay when using the TSM compared with that of the original schedule and MDM solution.}
	
	To study the quality of the solution over the entire parameter space, we vary one parameter in each \update{figure} that reports propagated delay comparisons. \update{Figure \ref{fig:budget}} reports a comparison for the \update{reschedule budget fractions in \{$0.25,0.5,0.75,1,2\}$}. Given a budget fraction, the corresponding reschedule budget is computed by multiplying the average value of the total primary flight delay of each of the $30$ recourse scenarios with the budget fraction value. \update{Figure \ref{fig:mean}} reports comparisons for distributions in \update{\{$Exponential(30)$, $LogNormal(30,15)$, $TruncatedNormal(30,15)$\}.} \update{Figure \ref{fig:mean}} fixes the distribution as exponential and reports comparisons for mean values of $\{15,30,45,60\}$ minutes. \update{Figure \ref{fig:scenarios} reports comparisons for the number of training scenarios in \{$10$,$20$,$30$,$40$,$50$\}. These figures show that the reduction in propagated delay achieved with TSM is significantly better than the original schedule and mean delay model, and that this reduction is agnostic of the underlying data.}

	\begin{table}
		\begin{center}
			\caption{Runtime comparison for column-generation strategies}
			\label{tab:columngen}
			\begin{tabular}{c c c c c} 
				\hline
				\textit{Instance}   &   \textit{Enumeration}  &  \textit{All paths} &  \textit{Best paths}  &  \textit{First Paths}    \\
				\hline
				s1   &   958.58   &  112.33  &   75.61  &  77.45    \\
				s2   &   161.19   &  63.45   &   47.46  &  49.87    \\
				s3   &   170.61   &  19.64   &   9.87   &  9.49     \\
				s4   &   417.46   &  28.32   &   15.2   &  14.28    \\
				s5   &   3086.92  &  121.61  &   65.81  &  69.34    \\
				s6   &   1860.32  &  1854.57 &   461.2  &  524.48   \\
				\hline
			\end{tabular}
		\end{center}
	\end{table}	
	
	In addition to the data-related parameters discussed so far, our approach also has several technical parameters, such as the type of column-generation strategy and the use of single versus multiple cuts for the $L$-shaped method. We use our final set of experiments to empirically select a set of these parameters that give the best runtime performance. The results are reported in Tables \ref{tab:columngen}, \ref{tab:threads}, \ref{tab:cut}, and \ref{tab:caching}.
	
	We obtain the values for each row in these tables as follows. First, we generate $30$ random delay scenarios using the default parameters specified in Section \ref{subsec:network}. Then we run Algorithm \ref{algo:l-shaped} for each value of the tested parameter and collect the solution time. We smooth out aberrations by repeating this $5$ times and reporting the average of these values as the time. The same procedure applies for values other than the solution time reported in Table \ref{tab:cut}.
	
	Table \ref{tab:columngen} reports a comparison between the different column-generation strategies presented in Section \ref{subsec:pricing-problem}. In this test, the first-paths and best-paths strategies are run with $N$=10, i.e., by selecting the first $10$ and the $10$ most negative reduced-cost columns, respectively. The results reported in this table are in line with the intuition that enumerating all columns should take much longer than using a delayed column-generation procedure with pricing. Among the pricing strategies, the best-paths and first-paths strategies are both clearly better than the all-paths strategy, which adds all negative reduced-cost columns to the restricted recourse problems.
	
	Table \ref{tab:threads} reports a run-time comparison with an increase in the number of threads. While it is indeed true that parallel solving should be faster, it is not practically obvious that this should be true. Specifically, we expected that the performance should stagnate or worsen when the number of threads exceeds the number of logical cores, but Table \ref{tab:threads} shows that this is not the case. Though the gain in performance declines with increasing threads, on an absolute basis, increasing the number of threads up to $30$ seems to improve the overall runtime. Increases beyond this are not helpful, as the maximum number of problems that can be solved in parallel is the number of recourse problems, which is $30$. Table \ref{tab:cut} reports a runtime comparison between the single- and multi-cut versions of Algorithm \ref{algo:l-shaped}. Clearly, the multi-cut version is better than the single-cut version in terms of the solution time, the Benders percentage gap (reported in the $Gap$ column), and the number of iterations. As the memory used to store and add cuts is minuscule in comparison to the rest of the data, the greater number of cuts in the multi-cut version does not affect performance at all. In Table \ref{tab:caching}, we present the results of caching the columns between the iterations for Algorithm \ref{algo:l-shaped}. We noticed that the columns generated in an iteration of the $L$-shaped method require only flight data and propagated delay data, and are unaffected by changes in the first-stage reschedule solution. This allows them to be cached and reused in future iterations, which in turn allows pricing problems to be warm-started with promising columns. As Table \ref{tab:caching} indicates, we were not able to find a clear advantage of this approach. While we certainly do not discard this idea, we recommend against using it, based purely on an ease-of-implementation perspective.

	\begin{table}
		\begin{center}
			\caption{Runtime comparison for multiple threads}
			\label{tab:threads}
			\begin{tabular}{c c c c c}
				\hline & \multicolumn{4}{c}{\textit{Number of parallel solvers}}\\
				\cline{2-5}
				\textit{Instance}   &   \textit{1}  &  \textit{10} &  \textit{20}  &  \textit{30}    \\
				\hline
				s1   &   692.71   &  123.49  &   98.63  &  77.88    \\
				s2   &   402.31   &  74.53   &   60.54  &  48.52    \\
				s3   &   64.64    &  12.58   &   10.16  &  8        \\
				s4   &   117.12   &  22.5    &   18.53  &  14.4     \\
				s5   &   607.55   &  104.76  &   88.55  &  74.04    \\
				s6   &   6490.08  &  986.12  &   736.96 &  519.8    \\
				\hline
			\end{tabular}
		\end{center}
	\end{table}

	\begin{table}
		\begin{center}
			\caption{Comparison of single- vs multi-cut L-shaped method}
			\label{tab:cut}
			\begin{tabular}{c c c c c | c c c c}
				\hline & \multicolumn{4}{c|}{\textit{Multi-cut}} &  \multicolumn{4}{c}{\textit{Single-cut}}\\
				\cline{2-9}
				\textit{Instance}   &   \textit{Time}  &  \textit{Gap} &  \textit{Cuts}  &  \textit{Iter}   &   \textit{Time}  &  \textit{Gap} &  \textit{Cuts}  &  \textit{Iter} \\
				\hline
				s1 & 686.81  & 0.4  & 883.4 & 30   & 708.91  & 24.28 & 30   & 30    \\
				s2 & 406.39  & 2.51 & 899.8 & 30   & 455.63  & 33.08 & 30   & 30    \\
				s3 & 58.16   & 0    & 85    & 4.8  & 214.2   & 0     & 19.4 & 19.6  \\
				s4 & 105.99  & 0.01 & 304.6 & 14   & 223.77  & 11.88 & 30   & 30    \\
				s5 & 579.02  & 0.03 & 340   & 14.4 & 1,172.27 & 12.36 & 30   & 30    \\
				s6 & 6,488.77 & 3.87 & 900   & 30   & 6,724.13 & 19.02 & 30   & 30    \\
				
				\hline
			\end{tabular}
		\end{center}
	\end{table}

	\begin{table}
		\begin{center}
			\caption{Runtime comparison for caching columns between iterations}
			\label{tab:caching}
			\begin{tabular}{c c c} 
				\hline
				\textit{Instance}   &   \textit{Caching}  &  \textit{No Caching} \\
				\hline
				s1   &   686.78  &  715.77    \\
				s2   &   399.28  &  422.96    \\
				s3   &   62.31   &  61.52     \\
				s4   &   112.8   &  105.6     \\
				s5   &   615.87  &  585.76    \\
				s6   &   6,372.9  &  6,198.92   \\
				\hline
			\end{tabular}
		\end{center}
	\end{table}	
	
	\vspace{10cm}
	\section{Conclusions and future research}\label{sec:conclusions}
	In this research, we present a two-stage stochastic programming model that adds time buffers to flight connections in order to make a schedule more robust to uncertain delays. By ``robust'', we mean that the schedule is more accommodating to changes in scheduled times and has fewer delays propagated to downstream flights. To solve the two-stage model, we present a solution framework that combines an outer approximation method with a delayed column-generation routine. We conduct a thorough qualitative and quantitative analysis of the proposed framework and report extensive computational results. To efficiently solve large-scale instances of the model, we adopt various software engineering techniques such as caching and parallelism. Our results highlight that the operational delay reduction can be significant using our proposed methodology compared to a deterministic approach.
	
	There are several interesting directions for extending this work, and we highlight a few \update{here}. First, the model can be made into a closer approximation of reality by considering more business constraints such as maintenance events and crew-friendliness. Another direction would be to study the scalability of our approach when more complex modifications such as cancellations, diversions, and overbooking are allowed in the first-stage. We have observed that, in practice, strategies to minimize delays can be quite diverse. While some airlines want to spread out delays among several flights to prohibit exorbitant delays for a single flight, other airlines want almost the exact opposite with the idea of minimizing the number of flights with delays. Making our model flexible enough to allow such variety in rescheduling and delay strategies is a worthwhile idea to pursue in the future. Also, from a modelling perspective, appropriate risk-averse objectives other than the risk-neutral expectation function can be evaluated in the second-stage.
	
	\section*{Acknowledgements}
	The authors would like to thank Sabre for providing the anonymized flight network data that we used for the computational studies in this article.
	
	\section*{Appendix}
	
	In this section, we report numbers used for the charts in Figures \ref{fig:budget}, \ref{fig:dist}, \ref{fig:mean} and \ref{fig:scenarios}.
	
	\update{Nomenclature common to all the following tables is listed below. All total propagated values are reported in minutes.
		\begin{itemize}
			\setlength\itemsep{1pt}
			\item \textit{Instance}: name of instance.
			\item \textit{Original}: average total propagated delay for the original schedules.
			\item \textit{MDM}: average total propagated delay with the schedule adjusted by the mean delay model solution.
			\item \textit{TSM}: average total propagated delay with the schedule adjusted by the TSM.
			\item \textit{RR over Original (\%)}: relative improvement of the MDM solution over the original ($100 \times (Original-TSM)/Original$).
			\item \textit{RR over MDM (\%)}: relative improvement of the TSM solution over the original ($100 \times (MDM-TSM)/MDM$). 
		\end{itemize}
	}	
	
	\begin{table}
		\begin{center}
			\caption{Total propagated delay improvements for different budgets}
			\label{tab:budget}
			\begin{tabular}{c c c c c c c}
				\hline & \multicolumn{6}{c}{\textit{Average total propagated delay}}\\
				\cline{3-7}
				\textit{Budget fraction}   &   \textit{Instance}  &  \textit{Original} &  \textit{MDM}  &  \textit{TSM}  &  \textit{RR over Original (\%)}  &  \textit{RR over MDM (\%)}    \\
				\hline
				0.25 & s1 & 845 & 628.06 & 562.57 & 33.42 & 10.43 \\
				& s2 & 850.82  & 611.65 & 520.17 & 38.86 & 14.96  \\
				& s3 & 50.24   & 26.88  & 15.68  & 68.79 & 41.67  \\
				& s4 & 219.37  & 145.93 & 135.86 & 38.07 & 6.9    \\
				& s5 & 254.29  & 215.02 & 160.18 & 37.01 & 25.5   \\
				& s6 & 1221.54 & 993.74 & 844.16 & 30.89 & 15.05  \\
				
				0.5 & s1 & 836.37 & 474.79 & 406.51 & 51.4 & 14.38  \\
				& s2 & 844.62  & 416.29 & 363.95 & 56.91 & 12.57 \\
				& s3 & 42.45   & 19.89  & 8.6    & 79.74 & 56.76 \\
				& s4 & 232.55  & 150.1  & 117.32 & 49.55 & 21.84 \\
				& s5 & 250.1   & 123.74 & 115.61 & 53.77 & 6.57  \\
				& s6 & 1231.86 & 799.37 & 672.06 & 45.44 & 15.93 \\
				
				0.75 & s1 & 861.65 & 373.57 & 365.71 & 57.56 & 2.1  \\
				& s2 & 868.94  & 345.26 & 303.68 & 65.05  & 12.04 \\
				& s3 & 46.81   & 25.88  & 11.76  & 74.88  & 54.56 \\
				& s4 & 218.15  & 132.55 & 87.93  & 59.69  & 33.66 \\
				& s5 & 242.06  & 116.37 & 102.03 & 57.85  & 12.32 \\
				& s6 & 1244.04 & 648.78 & 566.47 & 54.47  & 12.69 \\
				
				1 & s1 & 832.36 & 349.93 & 272.63 & 67.25 & 22.09  \\
				& s2 & 829.33  & 316.21 & 209.45 & 74.74  & 33.76 \\
				& s3 & 49.48   & 29.62  & 19.71  & 60.17  & 33.46 \\
				& s4 & 233.37  & 155.23 & 106.54 & 54.35  & 31.37 \\
				& s5 & 246.86  & 123.38 & 89.9   & 63.58  & 27.14 \\
				& s6 & 1197.72 & 505.05 & 502.65 & 58.03  & 0.48 \\
				
				2 & s1 & 849.18 & 351.68 & 238.15 & 71.96 & 32.28  \\
				& s2 & 851.63   & 344.38 & 222.88 & 73.83 & 35.28 \\
				& s3 & 49.12    &  28.81 & 16.94  & 65.51 & 41.2  \\
				& s4 & 222.53   & 144.08 & 95.3   & 57.17 & 33.86 \\
				& s5 & 243.47   & 116.92 & 79.63  & 67.29 & 31.89 \\
				& s6 & 1237.37  & 538.89 & 434.22 & 64.91 & 19.42 \\
				\hline
			\end{tabular}
		\end{center}
	\end{table}

	\begin{table}
		\begin{center}
			\caption{Total propagated delay improvements for different distributions}
			\label{tab:distribution}
			\begin{tabular}{c c c c c c c}
				\hline & \multicolumn{6}{c}{\textit{Average total propagated delay}}\\
				\cline{3-7}
				\textit{Distribution}   &   \textit{Instance}  &  \textit{Original} &  \textit{MDM}  &  \textit{TSM}  &  \textit{RR over Original (\%)}  &  \textit{RR over MDM (\%)}    \\
				\hline
				Exp(30) & s1 & 2,050.08 & 1,562.11 & 1,230.08 & 40 & 21.26 \\
				& s2 & 1,993.59  & 1336.85 & 1,107.84 & 44.43 & 17.13  \\
				& s3 & 141.43   & 87.52  & 55.89  & 60.48 & 36.14  \\
				& s4 & 701.25  & 434.87 & 391.68 & 44.15 & 9.93    \\
				& s5 & 599.99  & 411.45 & 330.68 & 44.89 & 19.63   \\
				& s6 & 3,281.99 & 2,542.58 & 2,113.81 & 35.59 & 16.86  \\
				
				LogNormal(30,15) & s1 & 1,966.24 & 1,233.31 & 867.31 & 55.89 & 29.68  \\
				& s2 & 1849.07  & 999.45 & 663.77 & 64.1 & 33.59 \\
				& s3 & 116.12   & 46.47  & 24.7    & 78.73 & 46.85 \\
				& s4 & 575.49  & 223.43  & 203.98 & 64.56 & 8.71 \\
				& s5 & 580.18   & 378.54 & 210.44 & 63.73 & 44.41  \\
				& s6 & 3,187.22 & 2,251.7 & 1619.3 & 49.19 & 28.09 \\
				
				TruncNormal(30,15) & s1 & 2,008.96 & 1,204.15 & 903.91 & 55.01 & 24.93  \\
				& s2 & 1919.41  & 900.75 & 693.00 & 63.84  & 22.95 \\
				& s3 & 115.87   & 39.72  & 18.16  & 84.33  & 54.28 \\
				& s4 & 615.21  & 238.11 & 207.77  & 66.23  & 16.26 \\
				& s5 & 580.18  & 378.54 & 210.44 & 63.73  & 44.41 \\
				& s6 & 3,187.22 & 2,251.7 & 1619.3 & 49.19  & 28.09 \\
				
				\hline
			\end{tabular}
		\end{center}
	\end{table}

	\begin{table}
		\begin{center}
			\caption{Total propagated delay improvements for different distribution means}
			\label{tab:mean}
			\begin{tabular}{c c c c c c c}
				\hline & \multicolumn{6}{c}{\textit{Average total propagated delay}}\\
				\cline{3-7}
				\textit{Distribution}   &   \textit{Instance}  &  \textit{Original} &  \textit{MDM}  &  \textit{TSM}  &  \textit{RR over Original (\%)}  &  \textit{RR over MDM (\%)}    \\
				\hline
				Exp(15) & s1 & 860.09 & 521.5 & 472.28 & 45.1 & 9.44 \\
				& s2 & 853.49  & 453.26 & 395.58 & 53.65 & 12.73  \\
				& s3 & 42.41   & 23.41  & 9.34  & 77.98 & 60.10  \\
				& s4 & 235.08  & 155.06 & 122.52 & 47.88 & 20.99    \\
				& s5 & 252.87  & 149.5 & 122.54 & 51.54 & 18.03   \\
				& s6 & 1,233.16 & 849.91 & 705.01 & 42.83 & 17.05  \\
				
				Exp(30) & s1 & 2,050.08 & 1,562.11 & 1,230.08 & 40 & 21.26 \\
				& s2 & 1,993.59 & 1,336.85 & 1,107.84 & 44.43 & 17.13 \\
				& s3 & 141.43 & 87.52 & 55.89 & 60.48 & 36.14 \\
				& s4 & 701.25 & 434.87 & 391.68 & 44.15 & 9.93 \\
				& s5 & 599.99 & 411.45 & 330.68 & 44.89 & 19.63 \\
				& s6 & 3,242.44 & 2,538.83 & 2,125.24 & 34.46 & 16.29 \\
				
				Exp(45) & s1 & 3,504.48 & 2,554.76 & 2,286.38 & 34.76 & 10.51 \\
				& s2 & 3,079.29 & 1,930.02 & 1,818.79 & 40.93 & 5.76 \\
				& s3 & 267.3 & 166.12 & 142.39 & 46.73 & 14.28 \\
				& s4 & 1,199.09 & 762.59 & 703.14 & 41.36 & 7.80 \\
				& s5 & 1,042.92 & 723.06 & 653.13 & 37.37 & 9.67 \\
				& s6 & 5,715.73 & 4,359.76 & 3846.91 & 32.7 & 11.76 \\
				
				Exp(60) & s1 & 5247.03 & 3922.66 & 3715.71 & 29.18 & 5.28 \\
				& s2 & 4674.16 & 3045.3 & 2938.05 & 37.14 & 3.52 \\
				& s3 & 412.07 & 280.51 & 257.87 & 37.42 & 8.07 \\
				& s4 & 1,825.04 & 1,322.64 & 1,168.95 & 35.95 & 11.62 \\
				& s5 & 1,437.66 & 1,138.58 & 958.5 & 33.33 & 15.82 \\
				& s6 & 8,822.83 &  6,750.83 &  6198.65 &  29.74 & 8.18 \\
				
				\hline
			\end{tabular}
		\end{center}
	\end{table}

	\begin{table}
		\begin{center}
			\caption{Total propagated delay improvements for different numbers of training scenarios}
			\label{tab:scenario}
			\begin{tabular}{c c c c c c c}
				\hline & \multicolumn{6}{c}{\textit{Average total propagated delay}}\\
				\cline{3-7}
				\textit{Scenarios}   &   \textit{Instance}  &  \textit{Original} &  \textit{MDM}  &  \textit{TSM}  &  \textit{RR over Original (\%)}  &  \textit{RR over MDM (\%)}    \\
				\hline
				
				10 & s1  & 1,960.17 &  1,115.96 &  846.1 &  56.84 &  24.18 \\
				20 &  &  1,945.01 &  1,146.31 &  830.28 &  57.31 &  27.57 \\
				30 &  &  1,941.91 &  1,112.12 &  831.43 &  57.18 &  25.24 \\
				40 &  &  1,943.8 &  1,123.35 &  812.97 &  58.18 &  27.63 \\
				50 &  &  1,949.21 &  1,117.64 &  806.86 &  58.61 &  27.81 \\
				
				10 &  s2  & 1,857.88 &  811.84 &  614.34 &  66.93 &  24.33 \\
				20 &  &  1,867.09 &  789.97 &  624.08 &  66.57 &  21 \\
				30 &  &  1,851.5 &  929.35 &  605.88 &  67.28 &  34.81 \\
				40 &  &  1,856.5 &  970.19 &  601.81 &  67.58 &  37.97 \\
				50 &  &  1,845.01 &  938.06 &  609.64 &  66.96 &  35.01 \\
				
				10  & s3  & 108.81 &  39.26 &  21.23 &  80.49 &  45.92 \\
				20  &  & 110.76 &  43.02 &  23.96 &  78.37 &  44.3 \\
				30  &  & 112.25 &  38.89 &  20.83 &  81.44 &  46.44 \\
				40 &  &  118.27 &  47.47 &  26.97 &  77.2 &  43.19 \\
				50 &   & 114.05 &  42.11 &  19.25 &  83.12 &  54.29 \\
				
				10 &  s4  & 587.27 &  212.83 &  183.25 &  68.8 &  13.9 \\
				20 &  &  583.26 &  189.41 &  161.75 &  72.27 &  14.6 \\
				30 &  &  610.97 &  202.41 &  184.06 &  69.87 &  9.07 \\
				40  &  & 580.37 &  194.53 &  163.15 &  71.89 &  16.13 \\
				50  &  & 586.65 &  188.69 &  168.23 &  71.32 &  10.84 \\
				
				10 &  s5  & 562.62 &  330.84 &  229.6 &  59.19 &  30.6 \\
				20 &  &  560.79 &  303.95 &  213.6 &  61.91 &  29.73 \\
				30  &  & 554.21 &  306.81 &  200.55 &  63.81 &  34.63 \\
				40 &  &  544.98 &  218.56 &  206.14 &  62.17 &  5.68 \\
				50 &  &  548.48 &  304.2 &  203.92 &  62.82 &  32.97 \\
				
				10 &  s6 &  3,053.12 &  2,143.16 &  1,508.29 &  50.6 &  29.62 \\
				20 &  &  3,082.73 &  2,221.58 &  1,560.49 &  49.38 &  29.76 \\
				30 &  &  3,062.17 &  2,137.81 &  1,527.92 &  50.1 &  28.53 \\
				40 &  &  3,099.17 &  2,169.34 &  1,518.15 &  51.01 &  30.02 \\
				50 &  &  3,078.29 &  2,086.7 &  1,487.11 &  51.69 &  28.73 \\
				
				\hline
			\end{tabular}
		\end{center}
	\end{table}	
	

	%
	%

	\bibliographystyle{spmpsci}      
	\bibliography{sample}   

\begin{thebibliography}{10}
\providecommand{\url}[1]{{#1}}
\providecommand{\urlprefix}{URL }
\expandafter\ifx\csname urlstyle\endcsname\relax
  \providecommand{\doi}[1]{DOI~\discretionary{}{}{}#1}\else
  \providecommand{\doi}{DOI~\discretionary{}{}{}\begingroup
  \urlstyle{rm}\Url}\fi

\bibitem{ahmadbeygi2010decreasing}
Ahmadbeygi, S., Cohn, A., Lapp, M.: Decreasing airline delay propagation by
  re-allocating scheduled slack.
\newblock IIE transactions \textbf{42}(7), 478--489 (2010)

\bibitem{arikan2013building}
Ar{\i}kan, M., Deshpande, V., Sohoni, M.: Building reliable air-travel
  infrastructure using empirical data and stochastic models of airline
  networks.
\newblock Operations Research \textbf{61}(1), 45--64 (2013)

\bibitem{balas1975set}
Balas, E., Padberg, M.: On the set-covering problem: Ii. an algorithm for set
  partitioning.
\newblock Operations Research \textbf{23}(1), 74--90 (1975)

\bibitem{balas1972set}
Balas, E., Padberg, M.W.: On the set-covering problem.
\newblock Operations Research \textbf{20}(6), 1152--1161 (1972)

\bibitem{bts}
BTS: {Bureau of Transportation Statistics}.
\newblock \url{https://www.transtats.bts.gov/OTDelay/OTDelayCause1.asp}

\bibitem{chiraphadhanakul2013robust}
Chiraphadhanakul, V., Barnhart, C.: Robust flight schedules through slack
  re-allocation.
\newblock EURO Journal on Transportation and Logistics \textbf{2}(4), 277--306
  (2013)

\bibitem{dunbar2012robust}
Dunbar, M., Froyland, G., Wu, C.L.: Robust airline schedule planning:
  Minimizing propagated delay in an integrated routing and crewing framework.
\newblock Transportation Science \textbf{46}(2), 204--216 (2012)

\bibitem{dunbar2014integrated}
Dunbar, M., Froyland, G., Wu, C.L.: An integrated scenario-based approach for
  robust aircraft routing, crew pairing and re-timing.
\newblock Computers \& Operations Research \textbf{45}, 68--86 (2014)

\bibitem{froyland2013recoverable}
Froyland, G., Maher, S.J., Wu, C.L.: The recoverable robust tail assignment
  problem.
\newblock Transportation Science \textbf{48}(3), 351--372 (2013)

\bibitem{kang2004degradable}
Kang, L.S.: Degradable airline scheduling: an approach to improve operational
  robustness and differentiate service quality.
\newblock Ph.D. thesis, Massachusetts Institute of Technology (2004)

\bibitem{klabjan2001robust}
Klabjan, D., Schaefer, A.J., Johnson, E.L., Kleywegt, A.J., Nemhauser, G.L.:
  Robust airline crew scheduling.
\newblock Proceedings of TRISTAN IV pp. 275--280 (2001)

\bibitem{lan2006planning}
Lan, S., Clarke, J.P., Barnhart, C.: Planning for robust airline operations:
  Optimizing aircraft routings and flight departure times to minimize passenger
  disruptions.
\newblock Transportation science \textbf{40}(1), 15--28 (2006)

\bibitem{marla2018robust}
Marla, L., Vaze, V., Barnhart, C.: Robust optimization: Lessons learned from
  aircraft routing.
\newblock Computers \& Operations Research \textbf{98}, 165--184 (2018)

\bibitem{rosenberger2004robust}
Rosenberger, J.M., Johnson, E.L., Nemhauser, G.L.: A robust fleet-assignment
  model with hub isolation and short cycles.
\newblock Transportation science \textbf{38}(3), 357--368 (2004)

\bibitem{shebalov2006robust}
Shebalov, S., Klabjan, D.: Robust airline crew pairing: Move-up crews.
\newblock Transportation science \textbf{40}(3), 300--312 (2006)

\bibitem{tahir2019integral}
Tahir, A., Desaulniers, G., El~Hallaoui, I.: Integral column generation for the
  set partitioning problem.
\newblock EURO Journal on Transportation and Logistics \textbf{8}(5), 713--744
  (2019)

\bibitem{talluri1996swapping}
Talluri, K.T.: Swapping applications in a daily airline fleet assignment.
\newblock Transportation Science \textbf{30}(3), 237--248 (1996)

\bibitem{van1969shaped}
Van~Slyke, R.M., Wets, R.: L-shaped linear programs with applications to
  optimal control and stochastic programming.
\newblock SIAM Journal on Applied Mathematics \textbf{17}(4), 638--663 (1969)

\bibitem{weide2010iterative}
Weide, O., Ryan, D., Ehrgott, M.: An iterative approach to robust and
  integrated aircraft routing and crew scheduling.
\newblock Computers \& Operations Research \textbf{37}(5), 833--844 (2010)

\bibitem{yan2016robust}
Yan, C., Kung, J.: Robust aircraft routing.
\newblock Transportation Science \textbf{52}(1), 118--133 (2016)

\bibitem{yen2006stochastic}
Yen, J.W., Birge, J.R.: A stochastic programming approach to the airline crew
  scheduling problem.
\newblock Transportation Science \textbf{40}(1), 3--14 (2006)

\end{thebibliography}


\begin{thebibliography}{}
	%
	%
	\bibitem{RefJ}
	Author, Article title, Journal, Volume, page numbers (year)
	\bibitem{RefB}
	Author, Book title, page numbers. Publisher, place (year)
	\end{thebibliography}

\end{document}